\RequirePackage{fix-cm}
\documentclass[smallextended, envcountsame]{svjour3}
\smartqed
\journalname{}

\usepackage[unicode = false,
    pdftoolbar = true,
    colorlinks = true,
    linkcolor = blue,
    citecolor = blue,
    filecolor = black,
    urlcolor = blue,
    breaklinks = true]{hyperref}

\usepackage[utf8]{inputenc}
\usepackage{listings}  % For Code 
\usepackage{enumitem}  
\usepackage{amsmath}   % For Math, duh
\usepackage{amssymb}   % For mathsfrs
\usepackage{xcolor}    % For Code Color
\usepackage{graphicx}  
\usepackage[ruled]{algorithm2e}
\usepackage{verbatim}
\usepackage{geometry}
\usepackage[justification=centering]{caption}
\usepackage{mathtools}
\usepackage{float}
\usepackage{bm}
\usepackage{subcaption}

\definecolor{gab}{HTML}{50c878}
\definecolor{pat}{HTML}{bf68ff}
\definecolor{sham}{HTML}{7faaff}

\DeclareMathAlphabet{\mymathbb}{U}{BOONDOX-ds}{m}{n}

\DeclareMathOperator*{\argmin}{argmin}

\makeatletter
\newcommand*{\rom}[1]{\expandafter\@slowromancap\romannumeral #1@}
\makeatother

\allowdisplaybreaks[4]

\title{%Floating point error analysis of general simplex derivatives\\
A general framework for floating point error analysis of simplex derivatives\\
}
\author{Yiwen Chen\and Warren Hare\and Amy Wiebe}
\institute{Department of Mathematics, University of British Columbia, Kelowna, British Columbia, V1V 1V7, Canada.\\
This research is partially funded by the Natural Sciences and Engineering Research Council of Canada (cette recherche est partiellement financ\'ee par le Conseil de recherches en sciences naturelles et en g\'enie du Canada), Discover Grant \#2023-03555.\\
\email{yiwchen@student.ubc.ca, warren.hare@ubc.ca, amy.wiebe@ubc.ca}}
\date{\today}

\begin{document}

    \maketitle

    \begin{abstract}
    Gradient approximations are a class of numerical approximation techniques that are of central importance in numerical optimization.  In derivative-free optimization, most of the gradient approximations, including the simplex gradient, centred simplex gradient, and adapted centred simplex gradient, are in the form of simplex derivatives.  Owing to machine precision, the approximation accuracy of any numerical approximation technique is subject to the influence of floating point errors.  In this paper, we provide a general framework for floating point error analysis of simplex derivatives.  Our framework is independent of the choice of the simplex derivative as long as it satisfies a general form.   We review the definition and approximation accuracy of the generalized simplex gradient and generalized centred simplex gradient.  We define and analyze the accuracy of a generalized version of the adapted centred simplex gradient.  As examples, we apply our framework to the generalized simplex gradient, generalized centred simplex gradient, and generalized adapted centred simplex gradient.  Based on the results, we give suggestions on the minimal choice of approximate diameter of the sample set.
    \end{abstract}
    \medskip


    \section{Introduction}\label{sec:intro}
    Derivative-free optimization (DFO) is a subfield of mathematical optimization that studies optimization algorithms that do not use derivatives of the objective or constraint functions \cite{audet2017derivative,conn2009introduction}.  DFO methods have been demonstrated to be effective and efficient in solving optimization problems where the derivatives are expensive or unavailable.  Two major categories of DFO methods are model-based methods and direct search methods.  In model-based methods, the objective and constraint functions are approximated by model functions that are constructed via derivative approximations.  Direct-search methods do not build any model functions and, instead, sample the objective function at a finite number of points in each iteration to search for a decrease.
    
    Gradient approximation techniques are of central importance in DFO methods.  Most of the gradient approximation techniques used in DFO methods are in the form of (first-order) simplex derivatives.  They can be computed by multiplying the inverse or pseudo-inverse of a matrix determined by sample points with a vector determined by function values on sample points (see the beginning of Section \ref{sec:fpe} for details).  For example, a linear approximation model in $\mathbb{R}^n$ can be constructed using $n+1$ well-poised sample points and the objective function values on these points~\cite{conn2009introduction}.  The gradient of this linear model is called the {\em simplex gradient} \cite[Chapter 9]{audet2017derivative}, which is an approximation of the true gradient of the objective function.  The approximation accuracy of the simplex gradient is $\mathcal{O}(\Delta)$, where $\Delta$ is the approximate diameter of the sample set~\cite{kelley1999iterative}.  Simplex derivatives are used in many more advanced model-based DFO methods including trust-region methods \cite[Chapter 11]{audet2017derivative}\cite{liuzzi2019trust} and quasi-Newton methods \cite{berahas2019derivative}.  Even in direct-search methods where, in general, no model functions are needed, simplex derivatives have been applied to improve convergence \cite{custodio2008using,kelley1999detection,kelley1999iterative}.
    
    Beyond the simplex gradient, various types of simplex derivatives with different approximation accuracy exist.  The {\em centred simplex gradient} \cite{hare2022error,kelley1999iterative} is constructed on two sets of points that are poised for linear interpolation and can be viewed as an average of two simplex gradients.  The major difference between the centred simplex gradient and the simplex gradient is that the centred simplex gradient is constructed on the original set of $n+1$ well-poised sample points and its reflection set through the reference point, hence uses a total number of $2n+1$ sample points.  This symmetry in the sample sets increases the approximation accuracy to $\mathcal{O}(\Delta^2)$~\cite{kelley1999iterative}.  Chen and Hare \cite{chen2023adapting} study the situation where the two sets used in the centred simplex gradient are not exactly symmetric.  They define a new simplex derivative called the {\em adapted centred simplex gradient} and show that it has an accuracy of $\mathcal{O}(\Theta\Delta+\Delta^2)$, where $\Theta$ is the maximum rotation angle from being exactly symmetric.

    In the three simplex derivatives mentioned above, it is also possible to use a different number of sample points.  In \cite{conn2008geometry,custodio2008using,custodio2007using,regis2015calculus}, the authors define the {\em generalized simplex gradient (GSG)}, which covers the cases where either fewer or more than $n+1$ sample points are provided, and provide the corresponding error bounds.  Hare et al. \cite{hare2022error} define and study the {\em generalized centred simplex gradient (GCSG)}, which does not require exactly $2n+1$ points.  Similarly, the analysis and results of the adapted centred simplex gradient \cite{chen2023adapting} can be generalized.  One contribution of this paper is to extend the results in \cite{chen2023adapting}.  In particular, this paper defines the {\em generalized adapted centred simplex gradient (GACSG)} and analyzes its approximation accuracy.
    
    For any numerical approximation technique, a key consideration is how the floating point errors impact its approximation accuracy.  As modern computers can only store real numbers in a finite number of bits, the real numbers are, in general, stored in an approximated representation with the relative approximation error bounded by the machine precision (or machine epsilon)~\cite{higham2002accuracy}.  This floating point error (or round-off error) may cause the approximation accuracy of numerical approximation techniques to be different than expected from a purely mathematical perspective.  For example, consider the following matrix $A$ and its true inverse $A^{-1}$:
    \begin{equation*}
        A=10
        \begin{bmatrix}
            ~1~~ & -\cos{\theta}\\
            ~0~~ & \sin{\theta}
        \end{bmatrix}~~~\text{and}~~~
        A^{-1}=
        0.1 \begin{bmatrix}
            ~1~~ &\displaystyle \cos{\theta}\slash\sin{\theta}\\
            ~0~~ & \displaystyle 1\slash\sin{\theta}
        \end{bmatrix}.
    \end{equation*}
    We let $\theta\in\{10^{-1},10^{-4},10^{-7},10^{-10},10^{-13},10^{-16}\}$ and use Python 3 with library {\tt numpy 1.26.4} to compute the inverse of $A$.  Note that {\tt numpy 1.26.4} uses the double-precision floating-point format according to IEEE 754, so the machine precision is approximately $2.2e-16$.  We represent the error by the spectral norm of the difference between the true inverse $A^{-1}$ and computed inverse $\overline{A^{-1}}$.  The following table shows how this error changes as $\theta$ decreases.
    \begin{center}
        \begin{tabular}{|c|c|c|c|c|c|c|}
        \hline
        $\theta$ & $10^{-1}$ & $10^{-4}$ & $10^{-7}$ & $10^{-10}$ & $10^{-13}$ & $10^{-16}$ \\ \hline
        $\left\|A^{-1}-\overline{A^{-1}}\right\|$ & $2.2e-16$ & $2.5e-13$ & $1.2e-10$ & $1.7e-07$ & $0.0e+00$ & $3.5e-01$ \\ \hline
        \end{tabular}
    \end{center}
    
    In this paper, we examine how the floating point error influences simplex derivatives and provide error bounds between the true gradient and the simplex derivatives given by computer calculations.  We provide a general framework for floating point error analysis of simplex derivatives.  Our framework is independent of the choice of simplex derivatives as long as they satisfy a general form (see the beginning of Section \ref{sec:fpe} for details).  In particular, we show how our framework can be applied to the GSG, GCSG, and GACSG, and give suggestions on the minimal choices of $\Delta$.
    
    The remainder of this paper is organized as follows.  We end this section with an introduction to our notation.  Section \ref{sec:def&appeb} provides the definitions of the GSG, GCSG, and GACSG, and gives their approximation errors without considering floating point errors.  Section \ref{sec:fpe} studies the floating point error that comes from computing the pseudo-inverse and evaluating the objective function, and gives a general error estimate between the true and computed simplex derivative.  Section~\ref{sec:exps} applies the general error estimate from Section \ref{sec:fpe} to the GSG, GCSG, and GACSG.  It also gives an example showing that if we know the structure of the simplex derivative, then it may be possible to obtain tighter error bounds than from our general error estimate.  Section~\ref{sec:deltamin} gives suggestions on the minimal choice of sample set approximate diameter $\Delta$ based on our error analysis.  Section \ref{sec:concl} summarizes this work and suggests some future research topics.

    \subsection{Notation}
        Let $\mymathbb{0}$ and $\mymathbb{1}$ denote the vector of all zeros and the vector of all ones with appropriate dimensions.  Let $I_n$ be the identity matrix of order $n$. For a matrix $A$, we use $A=(a_{ij})$ to mean the $ij$-th element of $A$ is $a_{ij}$, and denote its column space by $\mathrm{col}(A)$.  For real symmetric matrices, the function $\lambda_{\max}(\cdot)$ gives the maximum eigenvalue.  We use $\|\cdot\|$ to denote the spectral norm for matrices and the Euclidean norm for vectors.  For a set convex closed set $S$, function $\mathrm{proj}_S(\cdot)$ gives the Euclidean projection on $S$. The function $\mathcal{O}(\cdot)$ describes the limiting behaviour when the argument tends toward zero.  The set $$B_{\bar{\Delta}}(x)=\left\{y:\left\|y-x\right\|\le\bar{\Delta}\right\}$$ is the closed ball centred at $x$ with radius $\bar{\Delta}$.

        Throughout this paper, we assume $f:\mathbb{R}^n\to\mathbb{R}$.  For $k\in\{0,1,2,\ldots\}$, we say that $f\in\mathcal{C}^{k}$ on $B_{\bar{\Delta}}(x)$ if all partial derivatives of $f$ up to order $k$ exist and are continuous on $B_{\bar{\Delta}}(x)$.  Moreover, we say that a function $f\in\mathcal{C}^{k+}$ on $B_{\bar{\Delta}}(x)$ (with constant $\nu$) if $f\in\mathcal{C}^{k}$ on $B_{\bar{\Delta}}(x)$ and all the $k$-th order partial derivatives of $f$ are Lipschitz continuous on $B_{\bar{\Delta}}(x)$ (with Lipschitz constant $\nu$).  
    
        For a sample set $\mathbb{Y}=\{y_0,y_0+d_1,...,y_0+d_z\}$, define the $L$-matrix of $\mathbb{Y}$ by  
        \begin{equation}\label{eq:matL}
            L=L\left(\mathbb{Y}\right)=\left[y_0+d_1-y_0\cdots y_0+d_z-y_0\right]=\left[d_1\cdots d_z\right]\in\mathbb{R}^{n\times z}
        \end{equation}
        and define the approximate diameter of $\mathbb{Y}$ by
        \begin{equation*}
            \overline{\mathrm{diam}}(\mathbb{Y})=\max\left\{\left\|y_i-y_0\right\|:y_i\in\mathbb{Y}\right\}.
        \end{equation*}

        Recall that for a non-square or non-invertible square matrix, a generalization of the matrix inverse is the Moore–Penrose pseudo-inverse \cite{penrose1955generalized}.
        \begin{definition}[Moore–Penrose pseudo-inverse]
            For a matrix $A\in\mathbb{R}^{n\times z}$, the {\em Moore–Penrose pseudo-inverse} of $A$, denoted by $A^\dagger$, is the unique matrix in $\mathbb{R}^{z\times n}$ satisfying the following four conditions:
            \begin{enumerate}
                \item $AA^\dagger A=A,$
                \item $A^\dagger AA^\dagger=A^\dagger,$
                \item $(AA^\dagger)^\top=AA^\dagger,$
                \item $(A^\dagger A)^\top=A^\dagger A.$
            \end{enumerate}
        \end{definition}
        For any matrix $A\in\mathbb{R}^{n\times z}$, $A^\dagger$ exists and is unique \cite{penrose1955generalized}.  In particular, if $A$ has full rank, then $A^\dagger$ can be expressed by the following simple formulae \cite{penrose1955generalized}.
        \begin{enumerate}
            \item If $A\in\mathbb{R}^{n\times z}$ has full column rank $z$, then $$A^\dagger=(A^\top A)^{-1}A^\top$$ is a left-inverse of $A$, i.e., $A^\dagger A=I_z$. 
            \item If $A\in\mathbb{R}^{n\times z}$ has full row rank $n$, then $$A^\dagger=A^\top(AA^\top)^{-1}$$ is a right-inverse of $A$, i.e., $AA^\dagger=I_n$. 
        \end{enumerate}

\section{Simplex derivatives}\label{sec:def&appeb}
    In this section, we first introduce three types of simplex derivatives, which are the generalized simplex gradient, the generalized centred simplex gradient, and the generalized adapted centred simplex gradient.  Then, we provide their error bounds without floating point errors.

    The generalized simplex gradient generalizes the simplex gradient by allowing the number of sample points to be not exactly $n+1$.  This definition can be found in multiple papers, e.g., \cite{conn2008geometry,custodio2008using,custodio2007using,regis2015calculus}.
    \begin{definition}[Generalized Simplex Gradient, GSG]\label{def:gsg} 
        Suppose that $\mathbb{Y}=\{y_0,y_0+d_1,\ldots,y_0+d_p\}\subseteq\mathbb{R}^n$.  The {\em generalized simplex gradient (GSG)} of $f$ over $\mathbb{Y}$, denoted by $\nabla_{S}f(\mathbb{Y})$, is given by $$\nabla_{S}f(\mathbb{Y})=\left(L^\top\right)^\dagger\delta_{S}^{f(\mathbb{Y})},$$ where $L=L(\mathbb{Y})$ and 
        $$\delta_{S}^{f(\mathbb{Y})}=
        \begin{bmatrix}
        f(y_0+d_1)-f(y_0)\\
        \vdots\\
        f(y_0+d_p)-f(y_0)
        \end{bmatrix}.$$
    \end{definition}

    Next, we introduce the generalized centred simplex gradient, a thorough study of which can be found in \cite{hare2022error}.  Similar to the centred simplex gradient, the sample set of the generalized centred simplex gradient consists of an original set and its reflection set through the reference point.  However, in this generalized setting, we do not require exactly $n+1$ points in the original set.
    \begin{definition}[Generalized Centred Simplex Gradient, GCSG]\label{def:gcsg} 
        Suppose that $\mathbb{Y}^+=\{y_0,y_0+d_1,\ldots,y_0+d_p\}\subseteq\mathbb{R}^n$ and $\mathbb{Y}^-=\{y_0,y_0-d_1,\ldots,y_0-d_p\}\subseteq\mathbb{R}^n$. Let $\mathbb{Y}=\mathbb{Y}^+\cup\mathbb{Y}^-$. The {\em generalized centred simplex gradient (GCSG)} of $f$ over $\mathbb{Y}$, denoted by $\nabla_{CS}f(\mathbb{Y})$, is given by $$\nabla_{CS}f(\mathbb{Y})=\left(\left(L^+-L^-\right)^\top\right)^\dagger\delta_{CS}^{f(\mathbb{Y})},$$ where $L^+=L(\mathbb{Y}^+)$ and $L^-=L(\mathbb{Y}^-)$ are the $L$-matrices of $\mathbb{Y}^+$ and $\mathbb{Y}^-$, respectively, and
        $$\delta_{CS}^{f(\mathbb{Y})}=
        \begin{bmatrix}
        f(y_0+d_1)-f(y_0-d_1)\\
        \vdots\\
        f(y_0+d_p)-f(y_0-d_p)
        \end{bmatrix}.$$
    \end{definition}
    \begin{remark}
        By definition, we have $L^-=-L^+$.  Thus, an equivalent definition of the GCSG is $$\nabla_{CS}f(\mathbb{Y})=\left(\left(2L^+\right)^\top\right)^\dagger\delta_{CS}^{f(\mathbb{Y})}=\frac{1}{2}\left(\left(L^+\right)^\top\right)^\dagger\delta_{CS}^{f(\mathbb{Y})}.$$
    \end{remark}
    
    Notice that in the GCSG, the set $\mathbb{Y}^-$ is the exact reflection set of $\mathbb{Y}^+$ through $y_0$.  When $\mathbb{Y}^-$ is not the exact reflection set of $\mathbb{Y}^+$, a similar definition exists and is called the generalized adapted centred simplex gradient.  We denote the inexact reflection set of $\mathbb{Y}^+$ by $\widetilde{\mathbb{Y}}$.  The generalized adapted centred is defined based on the stretching and rotation relations between $\mathbb{Y}^+$ and $\widetilde{\mathbb{Y}}$.  The following two definitions are from \cite{chen2023adapting}.
    \begin{definition}[Stretching parameter]
        Suppose that $\mathbb{Y}^+=\{y_0,y_0+d_1,\ldots,y_0+d_p\}\subseteq\mathbb{R}^n$ and $\widetilde{\mathbb{Y}}=\{y_0,y_0-\widetilde{d}_1,\ldots,y_0-\widetilde{d}_p\}\subseteq\mathbb{R}^n$.  For all $i\in\{1,2,\ldots,p\}$, the {\em stretching parameter} $k_i$ is given by $$k_i=\frac{\left\|\widetilde{d}_i\right\|}{\left\|d_i\right\|}.$$
    \end{definition}
    \begin{definition}[Rotation angle and rotation matrix]
        Suppose that $\mathbb{Y}^+=\{y_0,y_0+d_1,\ldots,y_0+d_p\}\subseteq\mathbb{R}^n$ and $\widetilde{\mathbb{Y}}=\{y_0,y_0-\widetilde{d}_1,\ldots,y_0-\widetilde{d}_p\}\subseteq\mathbb{R}^n$.  For all $i\in\{1,2,\ldots,p\}$, the {\em rotation angle} $\theta_i$ is the angle between $d_i$ and $\widetilde{d}_i$, given by 
            $$\theta_i=\cos^{-1}\left(\frac{d_i^\top\widetilde{d}_i}{\left\|d_i\right\|\left\|\widetilde{d}_i\right\|}\right) \in [0, \pi].$$
        The {\em rotation matrix} $A_i$ is the unique matrix which rotates $\frac{d_i}{\left\|d_i\right\|}$ to $\frac{\widetilde{d}_i}{\left\|\widetilde{d}_i\right\|}$ (i.e., $\frac{\widetilde{d}_i}{\left\|\widetilde{d}_i\right\|}=A_i\frac{d_i}{\left\|d_i\right\|}$), and is in the form $A_i=P_i^\top A_i^\prime P_i$, where $P_i$ is an orthogonal matrix and
            $$A_i^\prime = 
            \begin{bmatrix}
        \cos(\theta_i) & -\sin(\theta_i) & \\
        \sin(\theta_i) & \cos(\theta_i)  & \\
        & & I_{n-2}
        \end{bmatrix}.$$
    \end{definition}
    
    Now we define the generalized adapted centred simplex gradient.  This is a generalization of the adapted centred simplex gradient defined and studied in \cite{chen2023adapting}.
    \begin{definition}[Generalized Adapted Centred Simplex Gradient, GACSG]\label{def:gacsg} 
        Suppose that $\mathbb{Y}^+=\{y_0,y_0+d_1,\ldots,y_0+d_p\}\subseteq\mathbb{R}^n$ and $\widetilde{\mathbb{Y}}=\{y_0,y_0-\widetilde{d}_1,\ldots,y_0-\widetilde{d}_p\}\subseteq\mathbb{R}^n$.  Let $\mathbb{Y}=\mathbb{Y}^+\cup\widetilde{\mathbb{Y}}$. The {\em generalized adapted centred simplex gradient (GACSG)} of $f$ over $\mathbb{Y}$, denoted by $\nabla_{ACS}f(\mathbb{Y})$, is given by $$\nabla_{ACS}f(\mathbb{Y})=\left(\left(L^+D-\widetilde{L}\right)^\top\right)^\dagger\delta_{ACS}^{f(\mathbb{Y})},$$ where $L^+=L(\mathbb{Y}^+)$ and $\widetilde{L}=L(\widetilde{\mathbb{Y}})$ are the $L$-matrices of $\mathbb{Y}^+$ and $\widetilde{\mathbb{Y}}$, respectively, 
        $$D=
        \begin{bmatrix}
        k_1^2 & & \\
         & \ddots & \\
        & & k_p^2
        \end{bmatrix},\ \mbox{and}\ 
        \delta_{ACS}^{f(\mathbb{Y})}=
        \begin{bmatrix}
        k_1^2\left(f(y_0+d_1)-f(y_0)\right)-\left(f(y_0-\widetilde{d}_1)-f(y_0)\right)\\
        \vdots\\
        k_p^2\left(f(y_0+d_p)-f(y_0)\right)-\left(f(y_0-\widetilde{d}_p)-f(y_0)\right)
        \end{bmatrix}.$$
    \end{definition}

    \subsection{Error bounds without floating point errors}
        Now we introduce the error bounds of the GSG, GCSG, and GACSG.  We note that all error bounds in this subsection do not consider the floating point error.  Error bounds with floating point errors are established in later sections.
        
        The error bounds of the GSG can be found in, e.g., \cite{conn2008geometry,custodio2008using,custodio2007using,regis2015calculus}.  Essentially, the error bound shows that the GSG provides an accuracy of $\mathcal{O}(\Delta)$ where $\Delta$ is the approximate diameter of the sample set.
        \begin{theorem}[Error bound of the GSG]\label{thm:GSGeb}
        	Suppose that $\mathbb{Y}=\{y_0,y_0+d_1,\ldots,y_0+d_p\}\subseteq\mathbb{R}^n$.  Let $\bar{\Delta}>0$.  Suppose that $f\in\mathcal{C}^{1+}$ on $B_{\bar{\Delta}}(y_0)$ with constant $\nu$ and $L\in\mathbb{R}^{n\times p}$ has full rank.  Suppose that $\Delta=\overline{\mathrm{diam}}(\mathbb{Y})\le\bar{\Delta}$.  Then
            	\begin{equation*}
            	    \left\|\mathrm{proj}_{\mathrm{col}(L)}(\nabla f(y_0))-\mathrm{proj}_{\mathrm{col}(L)}(\nabla_{S}f(\mathbb{Y}))\right\| \le \frac{\nu\sqrt{p}}{2}\left\|\widehat{L}^\dagger\right\|\Delta,
            	\end{equation*}
             where $\widehat{L}=L\slash\Delta$.  
             
             In particular, 
             \begin{equation*}
                 \mathrm{proj}_{\mathrm{col}(L)}(\nabla_Sf(\mathbb{Y})) = \mathrm{proj}_{\mathrm{col}(L)}(\nabla f(y_0)) + \mathcal{O}(\Delta).
             \end{equation*}
        \end{theorem}
        
        We benefit from the symmetry in the sample sets of the GCSG, which leads to an accuracy of $\mathcal{O}(\Delta^2)$.  The proof of the following error bound can be found in \cite[Theorems 3.3 and 3.6]{hare2022error}.  Notice that the proof only requires $f\in\mathcal{C}^{2+}$ which is a weaker assumption than the $f\in\mathcal{C}^{3}$ therein.
        \begin{theorem}[Error bound of the GCSG]\label{thm:GCSGeb}
        	Suppose that $\mathbb{Y}^+=\{y_0,y_0+d_1,\ldots,y_0+d_p\}\subseteq\mathbb{R}^n$ and $\mathbb{Y}^-=\{y_0,y_0-d_1,\ldots,y_0-d_p\}\subseteq\mathbb{R}^n$. Let $\mathbb{Y}=\mathbb{Y}^+\cup\mathbb{Y}^-$.  Let $\bar{\Delta}>0$.  Suppose that $f\in\mathcal{C}^{2+}$ on $B_{\bar{\Delta}}(y_0)$ with constant $\nu$ and $L^+\in\mathbb{R}^{n\times p}$ has full rank.  Suppose that $\Delta=\overline{\mathrm{diam}}(\mathbb{Y}^+)\le\bar{\Delta}$.  Then
            	\begin{equation*}
            	    \left\|\mathrm{proj}_{\mathrm{col}(L^+)}(\nabla f(y_0))-\mathrm{proj}_{\mathrm{col}(L^+)}(\nabla_{CS}f(\mathbb{Y}))\right\| \le \frac{\nu\sqrt{p}}{6}\left\|\left(\widehat{L}^+\right)^\dagger\right\|\Delta^2,
            	\end{equation*}
             where $\widehat{L}^+=L^+\slash\Delta$.
             
             In particular, 
             \begin{equation*}
                 \mathrm{proj}_{\mathrm{col}(L^+)}(\nabla_{CS}f(\mathbb{Y})) = \mathrm{proj}_{\mathrm{col}(L^+)}(\nabla f(y_0)) + \mathcal{O}(\Delta^2).
             \end{equation*}
        \end{theorem}

        For the GACSG, we note that all results in~\cite{chen2023adapting} only consider the determined case, i.e., the case where $L^+D-\widetilde{L}$ is invertible.  Now we generalize the results in \cite{chen2023adapting} to the case where $L^+D-\widetilde{L}$ only has full rank.  Following \cite{chen2023adapting}, we first introduce the following lemma.
        \begin{lemma}\label{lem:tle}
    	    Let $f\in\mathcal{C}^2$ and $x,d\in\mathbb{R}^n$. Then $$f(x+d)=f(x)+\nabla f(x)^\top d+\frac{1}{2}d^\top\nabla^2f(x)d+\int_0^1\int_0^1t_1d^\top\left(\nabla^2f(x+t_1t_2d)-\nabla^2f(x)\right)d\ \mathrm{d}t_2\mathrm{d}t_1.$$
    	\end{lemma}
    	\begin{proof}
    	    Using the fundamental theorem of calculus in $\mathbb{R}^n$ and the property of integrals, we have
    	    \begin{align*}
    	    &f(x+d)-f(x)-\nabla f(x)^\top d-\frac{1}{2}d^\top\nabla^2f(x)d \\
    	    =& \int_0^1\left(\nabla f(x+t_1d)-\nabla f(x)\right)^\top d\ \mathrm{d}t_1-\frac{1}{2}d^\top\nabla^2f(x)d\\
    	    =& \int_0^1\int_0^1t_1d^\top\nabla^2f(x+t_1t_2d)d\ \mathrm{d}t_2\mathrm{d}t_1-\frac{1}{2}d^\top\nabla^2f(x)d\\
    	    =& \int_0^1\int_0^1t_1d^\top\left(\nabla^2f(x+t_1t_2d)-\nabla^2f(x)\right)d\ \mathrm{d}t_2\mathrm{d}t_1.
    	    \end{align*}
    	    
    	    $\hfill\qed$
    	\end{proof}
     
        The proof of the following theorem is based on \cite[Theorems 3.3, 3.5, and 3.6]{chen2023adapting}.  We note that all matrix norms in \cite{chen2023adapting} are the Frobenius norm, while all matrix norms in this paper are the spectral norm.  However, the results in \cite{chen2023adapting} can be easily adapted to align with this paper since the Frobenius norm $\|\cdot\|_F$ and the spectral norm $\|\cdot\|$ have the relation $\|\cdot\|\le\|\cdot\|_F$.
        \begin{theorem}[Error bound of the GACSG]\label{thm:GACSGeb}
        	Suppose that $\mathbb{Y}^+=\{y_0,y_0+d_1,\ldots,y_0+d_p\}\subseteq\mathbb{R}^n$ and $\widetilde{\mathbb{Y}}=\{y_0,y_0-\widetilde{d}_1,\ldots,y_0-\widetilde{d}_p\}\subseteq\mathbb{R}^n$.  Let $\mathbb{Y}=\mathbb{Y}^+\cup\mathbb{Y}^-$.  Let $\bar{\Delta}>0$.  Suppose that $f\in\mathcal{C}^{2+}$ on $B_{\bar{\Delta}}(y_0)$ with constant $\nu$ and $(L^+D-\widetilde{L})\in\mathbb{R}^{n\times p}$ has full rank.  Suppose that $\Delta= \overline{\mathrm{diam}}(\mathbb{Y}^+) \le \bar{\Delta}$ and $\overline{\mathrm{diam}}(\widetilde{\mathbb{Y}})\le\bar{\Delta}$.  Then 
            \begin{align*}
                &~~\left\|\mathrm{proj}_{\mathrm{col}(L^+D-\widetilde{L})}(\nabla f(y_0))-\mathrm{proj}_{\mathrm{col}(L^+D-\widetilde{L})}(\nabla_{ACS}f(\mathbb{Y}))\right\|\\
                &\le \frac{1}{2}\max\left\{k_i^2\right\}\sqrt{p}\left(K_1\max\left\{\sin\theta_i\right\} + K_2\max\left\{\sqrt{1-\cos\theta_i}\right\}\right)\left\|\left(\widehat{L}^+D-\widehat{\widetilde{L}}\right)^\dagger\right\|\Delta\\
                &~~~+ \frac{\nu}{6}\max\left\{k_i^2\left(1+k_i\right)\right\}\sqrt{p}\left\|\left(\widehat{L}^+D-\widehat{\widetilde{L}}\right)^\dagger\right\|\Delta^2,
            \end{align*}
            where $$\widehat{L}^+=\frac{1}{\overline{\mathrm{diam}}(\mathbb{Y}^+)}L^+,\ \ \ \widehat{\widetilde{L}}=\frac{1}{\overline{\mathrm{diam}}(\mathbb{Y}^+)}\widetilde{L},$$
            $$K_1=\max\left\{\sqrt{\left(h_{11}^i-h_{22}^i\right)^2+4\left(h_{12}^i\right)^2}\right\},$$ and 
            $$K_2=\max\left\{\sqrt{\left(\sum\limits_{l=3}^n\left(h_{1l}^i\right)^2+\sum\limits_{l=3}^n\left(h_{2l}^i\right)^2\right)+\sqrt{\left(\sum\limits_{l=3}^n\left(h_{1l}^i\right)^2-\sum\limits_{l=3}^n\left(h_{2l}^i\right)^2\right)^2+4\left(\sum\limits_{l=3}^nh_{1l}^ih_{2l}^i\right)^2}}\right\},$$ where $h_{kl}^i$ is the $(kl)$-th element of matrix $P_i\nabla^2f(y_0)P_i^\top$.
            
            In particular, if all $k_i$ are near 1, then we have 
            \begin{equation*}
                \mathrm{proj}_{\mathrm{col}(L^+D-\widetilde{L})}(\nabla_{ACS}f(\mathbb{Y})) = \mathrm{proj}_{\mathrm{col}(L^+D-\widetilde{L})}(\nabla f(y_0)) + \mathcal{O}(\Theta \Delta + \Delta^2)
            \end{equation*}
            where $\Theta = \max\{\theta_i\}$.
        \end{theorem}
        \begin{proof}
            Since $L^+D-\widetilde{L}$ has full rank, we break our analysis into two cases: the case where $L^+D-\widetilde{L}$ has full-row rank and the case where $L^+D-\widetilde{L}$ has full-column rank.
            
            {\bf Case \rom{1}:} Suppose that $L^+D-\widetilde{L}$ has full-row rank.\\
            Select $i\in\{1,\ldots,p\}$.  Let $e_i$ be the $i$-th coordinate vector in $\mathbb{R}^p$.  Following the proof of \cite[Theorems~3.3, 3.5, and 3.6]{chen2023adapting}, we have
            \begin{align*}
                &~~\left\|e_i^\top\left(L^+D-\widetilde{L}\right)^\top\nabla f(y_0) - e_i^\top\delta_{ACS}^{f(\mathbb{Y})}\right\|\\
                &= \left\|\left(k_i^2d_i+k_iA_id_i\right)^\top\nabla f(y_0) - k_i^2\left(f(y_0+d_i)-f(y_0)\right)+\left(f(y_0-\widetilde{d}_i)-f(y_0)\right)\right\|\\
                &\le \frac{1}{2}k_i^2\left(\kappa_1^i\sin\theta_i + \kappa_2^i\sqrt{1-\cos\theta_i}\right)\Delta^2+\frac{\nu}{6}k_i^2\left(1+k_i\right)\Delta^3,
            \end{align*}
            where $$\kappa_1^i= \sqrt{\left(h_{11}^i-h_{22}^i\right)^2+4\left(h_{12}^i\right)^2},$$ and $$\kappa_2^i = \sqrt{\left(\sum\limits_{l=3}^n\left(h_{1l}^i\right)^2+\sum\limits_{l=3}^n\left(h_{2l}^i\right)^2\right)+\sqrt{\left(\sum\limits_{l=3}^n\left(h_{1l}^i\right)^2-\sum\limits_{l=3}^n\left(h_{2l}^i\right)^2\right)^2+4\left(\sum\limits_{l=3}^nh_{1l}^ih_{2l}^i\right)^2}}.$$
            Therefore,
            \begin{align*}
                &~~\left\|\left(L^+D-\widetilde{L}\right)^\top\nabla f(y_0) - \delta_{ACS}^{f(\mathbb{Y})}\right\|\\
                &\le \left\|\left(L^+D-\widetilde{L}\right)^\top\nabla f(y_0) - \delta_{ACS}^{f(\mathbb{Y})}\right\|_F\\
                &= \sqrt{\sum\limits_{i=1}^p\left\|e_i^\top\left(L^+D-\widetilde{L}\right)^\top\nabla f(y_0) - e_i^\top\delta_{ACS}^{f(\mathbb{Y})}\right\|^2}\\
                &\le \frac{1}{2}\max\left\{k_i^2\right\}\sqrt{p}\left(K_1\max\left\{\sin\theta_i\right\} + K_2\max\left\{\sqrt{1-\cos\theta_i}\right\}\right)\Delta^2+\frac{\nu}{6}\max\left\{k_i^2\left(1+k_i\right)\right\}\sqrt{p}\Delta^3,
            \end{align*}
            where $\|\cdot\|_F$ is the Frobenius norm.
            
            Notice that $((L^+D-\widetilde{L})^\top)^\dagger(L^+D-\widetilde{L})^\top=I_{n}$ and $\mathrm{col}(L^+D-\widetilde{L})=\mathbb{R}^n$.  We obtain
            \begin{align*}
                &~~\left\|\mathrm{proj}_{\mathrm{col}(L^+D-\widetilde{L})}(\nabla f(y_0))-\mathrm{proj}_{\mathrm{col}(L^+D-\widetilde{L})}(\nabla_{ACS}f(\mathbb{Y}))\right\|\\
                &= \left\|\nabla f(y_0)-\nabla_{ACS}f(\mathbb{Y})\right\|\\
                &= \left\|\left(\left(L^+D-\widetilde{L}\right)^\top\right)^\dagger\left(L^+D-\widetilde{L}\right)^\top\nabla f(y_0) - \left(\left(L^+D-\widetilde{L}\right)^\top\right)^\dagger\delta_{ACS}^{f(\mathbb{Y})}\right\|\\
                &\le \left\|\left(L^+D-\widetilde{L}\right)^\dagger\right\|\left\|\left(L^+D-\widetilde{L}\right)^\top\nabla f(y_0) - \delta_{ACS}^{f(\mathbb{Y})}\right\|\\
                &= \frac{1}{\Delta}\left\|\left(\widehat{L}^+D-\widehat{\widetilde{L}}\right)^\dagger\right\|\left\|\left(L^+D-\widetilde{L}\right)^\top\nabla f(y_0) - \delta_{ACS}^{f(\mathbb{Y})}\right\|\\
                &\le \frac{1}{2}\max\left\{k_i^2\right\}\sqrt{p}\left(K_1\max\left\{\sin\theta_i\right\} + K_2\max\left\{\sqrt{1-\cos\theta_i}\right\}\right)\left\|\left(\widehat{L}^+D-\widehat{\widetilde{L}}\right)^\dagger\right\|\Delta\\
                &~~~+ \frac{\nu}{6}\max\left\{k_i^2\left(1+k_i\right)\right\}\sqrt{p}\left\|\left(\widehat{L}^+D-\widehat{\widetilde{L}}\right)^\dagger\right\|\Delta^2.
            \end{align*}

            {\bf Case \rom{2}:} Suppose that $L^+D-\widetilde{L}$ has full-column rank.\\
            Notice that $(L^+D-\widetilde{L})^\top((L^+D-\widetilde{L})^\top)^\dagger=I_{p}$.  Using a similar proof as \cite[Theorems~3.3, 3.5, and~3.6]{chen2023adapting}, we get
            \begin{align*}
                \left\|\left(L^+D-\widetilde{L}\right)^\top\left(\nabla f(y_0) - \nabla_{ACS}f(\mathbb{Y})\right)\right\| \le& \frac{1}{2}\max\left\{k_i^2\right\}\sqrt{p}\left(K_1\max\left\{\sin\theta_i\right\} + K_2\max\left\{\sqrt{1-\cos\theta_i}\right\}\right)\Delta^2\\
                &+ \frac{\nu}{6}\max\left\{k_i^2\left(1+k_i\right)\right\}\sqrt{p}\Delta^3.
            \end{align*}
            Therefore,
            \begin{align*}
                &~~\left\|\mathrm{proj}_{\mathrm{col}(L^+D-\widetilde{L})}(\nabla f(y_0))-\mathrm{proj}_{\mathrm{col}(L^+D-\widetilde{L})}(\nabla_{ACS}f(\mathbb{Y}))\right\|\\
                &= \left\|\left(\left(L^+D-\widetilde{L}\right)^\top\right)^\dagger\left(L^+D-\widetilde{L}\right)^\top\left(\nabla f(y_0)-\nabla_{ACS}f(\mathbb{Y})\right)\right\|\\
                &\le \left\|\left(L^+D-\widetilde{L}\right)^\dagger\right\|\left\|\left(L^+D-\widetilde{L}\right)^\top\left(\nabla f(y_0)-\nabla_{ACS}f(\mathbb{Y})\right)\right\|\\
                &= \frac{1}{\Delta}\left\|\left(\widehat{L}^+D-\widehat{\widetilde{L}}\right)^\dagger\right\|\left\|\left(L^+D-\widetilde{L}\right)^\top\left(\nabla f(y_0)-\nabla_{ACS}f(\mathbb{Y})\right)\right\|\\
                &\le \frac{1}{2}\max\left\{k_i^2\right\}\sqrt{p}\left(K_1\max\left\{\sin\theta_i\right\} + K_2\max\left\{\sqrt{1-\cos\theta_i}\right\}\right)\left\|\left(\widehat{L}^+D-\widehat{\widetilde{L}}\right)^\dagger\right\|\Delta\\
                &~~~+ \frac{\nu}{6}\max\left\{k_i^2\left(1+k_i\right)\right\}\sqrt{p}\left\|\left(\widehat{L}^+D-\widehat{\widetilde{L}}\right)^\dagger\right\|\Delta^2.
            \end{align*}
         
        $\hfill\qed$
        \end{proof}

\section{Floating point errors}\label{sec:fpe}    
    Now we begin to analyze the floating point errors that may occur in gradient approximations.  Our goal is to establish error bounds of all types of simplex derivatives satisfying a general form.  Starting from this section, we label the points in $\mathbb{Y}$ by $\mathbb{Y}=\{y_0,\ldots,y_q\}$ and assume a general simplex derivative of the form:
    \begin{equation*}
        \nabla_Xf(\mathbb{Y})=(A^\top)^\dagger B\mymathbb{f}(\mathbb{Y}),
    \end{equation*}
    where $A\in\mathbb{R}^{n\times z}$, $B\in\mathbb{R}^{z\times (q+1)}$, and 
    \begin{equation*}
        \mymathbb{f}(\mathbb{Y})=\begin{bmatrix}
            f(y_0)\\
            \vdots\\
            f(y_q)
        \end{bmatrix}.
    \end{equation*}

    Notice that each of the GSG, GCSG, and GACSG has this form, where the necessary parameters are given in the following table.
    \begin{center}
    \begin{tabular}{cccccc}
        \hline
        Type & $\mathbb{Y}\subseteq\mathbb{R}^n$ & $z$ & $q$ & $A\in\mathbb{R}^{n\times z}$ & $B\in\mathbb{R}^{z\times (q+1)}$ \\ \hline
        GSG & $\{y_0,y_0+d_1,\ldots,y_0+d_p\}$ & $p$ & $p$ & $L$ & $[-\mymathbb{1}~I_p]$ \\
        GCSG & \begin{tabular}{c}$\{y_0,y_0+d_1,\ldots,y_0+d_p,$ \\ $\hspace{0.75cm}y_0-d_1,\ldots,y_0-d_p\}$\end{tabular} & $p$ & $2p$ & $2L^+$ & $[\mymathbb{0}~I_p~-I_p]$ \\
        GACSG & \begin{tabular}{c}$\{y_0,y_0+d_1,\ldots,y_0+d_p,$\\ $\hspace{0.75cm}y_0-\widetilde{d}_1,\ldots,y_0-\widetilde{d}_p\}$\end{tabular} & $p$ & $2p$ & $L^+D-\widetilde{L}$ &  $\begin{bmatrix}
            1-k_1^2 & k_1^2 & \cdots & 0 & -1 & \cdots & 0\\
            \vdots & \vdots & \ddots & \vdots & \vdots & \ddots & \vdots\\
            1-k_p^2 & 0 & \cdots & k_p^2 & 0 & \cdots & -1
        \end{bmatrix}$\\ \hline
    \end{tabular}        
    \end{center}

    Two major sources of floating point errors in computing simplex derivatives in the general form above are computing the pseudo-inverse and evaluating function values on the sample points.  In Subsection \ref{subsec:invfpe}, we analyze the floating point error in computing the pseudo-inverse.  In Subsection \ref{subsec:fevalfpe}, we analyze the floating point error in evaluating function values.  Subsection~\ref{subsec:invfevalfpe} combines both sources of floating point errors and gives an error bound between the true and computed simplex derivatives in the general form.

    \subsection{Floating point error in matrix inversion}\label{subsec:invfpe}
        The magnitude of floating point error in computing the pseudo-inverse is related to the numerical stability of the algorithms.  Our analysis is based on a class of algorithms that are sometimes referred to as numerically stable algorithms \cite{higham2002accuracy}.  Many classic algorithms are proven to satisfy this definition, e.g., the Householder bidiagonalization \cite{byers2008new} and the singular value decomposition~\cite{smoktunowicz2012numerical}.
        \begin{definition}
            Let $M\in\mathbb{R}^{n\times z}$.  An algorithm for computing $M^\dagger$ is {\em mixed forward-backward stable}, if the computed result $\overline{M^\dagger}\in\mathbb{R}^{z\times n}$ satisfies 
            \begin{equation*}
                \overline{M^\dagger} + \widetilde{M^\dagger} = \left(M + \widetilde{M}\right)^\dagger, 
            \end{equation*}
            for some $\widetilde{M^\dagger}\in\mathbb{R}^{z\times n}$ and $\widetilde{M}\in\mathbb{R}^{n\times z}$ with
            \begin{equation*}
                \left\|\widetilde{M}\right\|\le C\epsilon^*\left\|M\right\|~~\text{and}~~\left\|\widetilde{M^\dagger}\right\| \le C\epsilon^*\left\|\overline{M^\dagger}\right\|,
            \end{equation*}
            where $C$ is a small constant depending upon $n$ and $z$, and $\epsilon^*$ is machine precision.
        \end{definition}

        The following lemma provides an upper bound for the error between the true pseudo-inverse and the computed pseudo-inverse using a mixed forward-backward stable algorithm.
        \begin{lemma}\label{lem:invne}
            Let $\epsilon^*$ be the machine precision.  Suppose that $A\in\mathbb{R}^{n\times z}$ has full-column rank.  Suppose that $\overline{(A^\top)^\dagger}$ is computed by a mixed forward-backward stable algorithm with constant $C$. Let $\kappa(A)=\|A^\dagger\|\|A\|$ be the condition number of $A$.  If $C\epsilon^*\kappa(A)<1$, then
            \begin{align*}
                \left\|\left(A^\top\right)^\dagger - \overline{(A^\top)^\dagger}\right\| \le \frac{C\left\|A^\dagger\right\|\epsilon^*}{1-C\kappa(A)\epsilon^*}\left(\sqrt{2}\kappa(A)+\frac{1}{1-C\epsilon^*}\right).
            \end{align*}
        \end{lemma}
        \begin{proof}
            Notice that $\mathrm{rank}(A^\top)=z=\min\{n, z\}$ and $\|A^\dagger\|\|\widetilde{A^\top}\|\le C\epsilon^*\kappa(A)<1$.  From \cite[p.~29]{bjorck1996numerical} we have $\mathrm{rank}(A^\top)=\mathrm{rank}(A^\top+\widetilde{A^\top})$.  From here, we can complete the proof following \cite[Theorem~2.1]{smoktunowicz2012numerical}, noticing that the results therein hold for all matrices with full rank. 
    
        $\hfill\qed$
        \end{proof}
    
        Based on Lemma \ref{lem:invne}, the floating point error that comes from computing the pseudo-inverse in the simplex derivatives can be bounded as follows.
        \begin{theorem}[Floating point error in matrix inversion]\label{thm:invne}
            Let $\epsilon^*$ be the machine precision.  Suppose that $A\in\mathbb{R}^{n\times z}$ has full-column rank.  Suppose that $\overline{(A^\top)^\dagger}$ is computed by a mixed forward-backward stable algorithm with constant $C$. Let $\kappa(A)=\|A^\dagger\|\|A\|$ be the condition number of $A$.  Let $\overline{\nabla_Xf}(\mathbb{Y})=\overline{(A^\top)^\dagger}B\mymathbb{f}(\mathbb{Y})$.  If $C\epsilon^*\kappa(A)<1$, then
            \begin{equation*}
                \left\|\nabla_Xf(\mathbb{Y})-\overline{\nabla_Xf}(\mathbb{Y})\right\| \le \frac{C\left\|A^\dagger\right\|\epsilon^*}{1-C\kappa(A)\epsilon^*}\left(\sqrt{2}\kappa(A)+\frac{1}{1-C\epsilon^*}\right)\left\|B\right\|\left\|\mymathbb{f}(\mathbb{Y})\right\|.
            \end{equation*}
        \end{theorem}
        \begin{proof}
            Notice that 
            \begin{equation*}
                \left\|\nabla_Xf(\mathbb{Y})-\overline{\nabla_Xf}(\mathbb{Y})\right\| = \left\|\left(A^\top\right)^\dagger B\mymathbb{f}(\mathbb{Y}) - \overline{(A^\top)^\dagger}B\mymathbb{f}(\mathbb{Y})\right\|.
            \end{equation*}
            The proof is complete by using Lemma \ref{lem:invne} and the triangle inequality.
    
        $\hfill\qed$
        \end{proof}

    \subsection{Floating point error in function evaluations}\label{subsec:fevalfpe}
        In this subsection, we analyze the floating point error in objective function evaluations, which is the other source of floating point errors in computing simplex derivatives in the general form.  We note that the analysis in this section can also be viewed as analyzing the influence of noise in the objective function.
    
        For any $x$, we suppose that the computed function value of $f$ at $x$ has the form $\overline{f}(x)=(1+\epsilon_x)f(x)$, where the absolute value of $\epsilon_x$ is bounded by a constant.  In this paper, the constant is the machine precision $\epsilon^*$.  To adapt the analysis in this section to noisy functions, one should take the constant to be the noise level.  We note that some other research may define an alternative form $\overline{f}(x)=f(x)+\epsilon_x$.  Our analysis and results can be easily extended to that case.
    
        The following lemma examines the floating point error in computing the vector $\mymathbb{f}(\mathbb{Y})$.
        \begin{lemma}\label{lem:fevalne}
            Let $\epsilon^*$ be the machine precision.  For $x\in\mathbb{R}^n$, let $\overline{f}(x)=(1+\epsilon_x)f(x)$ with $|\epsilon_x|\le\epsilon^*$.  Define $f_M=\max\{|f(y)|:y\in\mathbb{Y}\}$ and 
            \begin{equation*}
                \overline{\mymathbb{f}}(\mathbb{Y})=\begin{bmatrix}
                    \overline{f}(y_0)\\
                    \vdots\\
                    \overline{f}(y_q)
                \end{bmatrix}.
            \end{equation*}
            Then
            \begin{equation*}
                \left\|\mymathbb{f}(\mathbb{Y})-\overline{\mymathbb{f}}(\mathbb{Y})\right\| \le \sqrt{q+1}f_M\epsilon^*.
            \end{equation*}
        \end{lemma}
        \begin{proof}
            Notice that $|\epsilon_{y_i}|\le\epsilon^*$ and $|f(y_i)|\le f_M$ for all $i=0,\ldots,q$.  We have
            \begin{equation*}
                \left\|\mymathbb{f}(\mathbb{Y})-\overline{\mymathbb{f}}(\mathbb{Y})\right\| = \sqrt{\sum\limits_{i=0}^q\left|f(y_i)-\overline{f}(y_i)\right|^2} = \sqrt{\sum\limits_{i=0}^q\left|\epsilon_{y_i}\right|^2\left|f(y_i)\right|^2} \le \sqrt{q+1}f_M\epsilon^*.
            \end{equation*}
    
        $\hfill\qed$
        \end{proof}
    
        Based on Lemma \ref{lem:fevalne}, the floating point error that comes from evaluating objective function values in the simplex derivatives can be bounded as follows.
        \begin{theorem}[Floating point error in function evaluations]\label{thm:fevalne}
            Let $\epsilon^*$ be the machine precision.  For $x\in\mathbb{R}^n$, let $\overline{f}(x)=(1+\epsilon_x)f(x)$ with $|\epsilon_x|\le\epsilon^*$.  Let $\overline{\nabla_X f}(\mathbb{Y})=(A^\top)^\dagger B\overline{\mymathbb{f}}(\mathbb{Y})$.  Define $f_M=\max\{|f(y)|:y\in\mathbb{Y}\}$. Then
            \begin{equation*}
                \left\|\nabla_X f(\mathbb{Y})-\overline{\nabla_X f}(\mathbb{Y})\right\| \le \sqrt{q+1}f_M\left\|A^\dagger\right\|\left\|B\right\|\epsilon^*.
            \end{equation*}
        \end{theorem}
        \begin{proof}
            The proof is complete by using Lemma \ref{lem:fevalne} and the triangle inequality.
    
        $\hfill\qed$
        \end{proof}

    \subsection{Floating point errors in matrix inversion and function evaluations}\label{subsec:invfevalfpe}
        In this subsection, we establish an error bound between the true and computed simplex derivatives in the general form, which combines both sources of errors analyzed in Subsections~\ref{subsec:invfpe} and~\ref{subsec:fevalfpe}. We then apply this error bound to three specific simplex derivatives, i.e., the GSG, GCSG, and GACSG, as examples.
    
        The following lemma gives an upper bound for the norm of the computed pseudo-inverse.  The proof is modified from the proof of \cite[Theorem 2.1]{smoktunowicz2012numerical}.
        \begin{lemma}\label{lem:bdpinv}
            Let $\epsilon^*$ be the machine precision.  Suppose that $A\in\mathbb{R}^{n\times z}$ has full-column rank.  Suppose that $\overline{(A^\top)^\dagger}$ is computed by a mixed forward-backward stable algorithm with constant $C$. Let $\kappa(A)=\|A^\dagger\|\|A\|$ be the condition number of $A$.  If $C\epsilon^*\kappa(A)<1$, then 
            \begin{equation*}
                \left\|\overline{(A^\top)^\dagger}\right\| \le \frac{\|A^\dagger\|}{(1-C\epsilon^*)(1-C\kappa(A)\epsilon^*)}.
            \end{equation*}
        \end{lemma}
        \begin{proof}
            By definition of the mixed forward-backward stable algorithm, we have
            \begin{equation*}
                \overline{(A^\top)^\dagger} = \left(A^\top+\widetilde{A^\top}\right)^\dagger - \widetilde{(A^\top)^\dagger},
            \end{equation*}
            with
            \begin{equation*}
                \left\|\widetilde{(A^\top)^\dagger}\right\|\le C\epsilon^*\left\|\overline{(A^\top)^\dagger}\right\|.
            \end{equation*}
            
            From \cite[Lemma 2.2, Equation (2.2)]{smoktunowicz2012numerical} with $\epsilon=C\epsilon^*$, we have
            \begin{equation*}
                \left\|\left(A^\top+\widetilde{A^\top}\right)^\dagger\right\| \le \frac{\|A^\dagger\|}{1-C\kappa(A)\epsilon^*}.
            \end{equation*}
            Therefore,
            \begin{equation*}
                \left\|\overline{(A^\top)^\dagger}\right\| \le \left\|\left(A^\top+\widetilde{A^\top}\right)^\dagger\right\| + \left\|\widetilde{(A^\top)^\dagger}\right\| \le \frac{\|A^\dagger\|}{1-C\kappa(A)\epsilon^*} + C\epsilon^*\|\overline{(A^\top)^\dagger}\|,
            \end{equation*}
            i.e.,
            \begin{equation*}
                \left(1-C\epsilon^*\right)\left\|\overline{(A^\top)^\dagger}\right\| \le \frac{\|A^\dagger\|}{1-C\kappa(A)\epsilon^*}.
            \end{equation*}
            Since $C\epsilon^*\kappa(A)<1$, we have
            \begin{equation*}
                \left\|\overline{(A^\top)^\dagger}\right\| \le \frac{\|A^\dagger\|}{(1-C\epsilon^*)(1-C\kappa(A)\epsilon^*)}.
            \end{equation*}
    
        $\hfill\qed$
        \end{proof}
    
        Now we combine the floating points errors from pseudo-inverse computation (Subsection \ref{subsec:invfpe}) and objective function evaluation (Subsection \ref{subsec:fevalfpe}) to get the following theorem.
        \begin{theorem}[Floating point errors in matrix inversion and function evaluations]\label{thm:invfevalne}
            Let $\epsilon^*$ be the machine precision.  For $x\in\mathbb{R}^n$, let $\overline{f}(x)=(1+\epsilon_x)f(x)$ with $|\epsilon_x|\le\epsilon^*$.  Suppose that $A\in\mathbb{R}^{n\times z}$ has full-column rank.  Suppose that $\overline{(A^\top)^\dagger}$ is computed by a mixed forward-backward stable algorithm with constant $C$. Let $\kappa(A)=\|A^\dagger\|\|A\|$ be the condition number of $A$.  Define $f_M=\max\{|f(y)|:y\in\mathbb{Y}\}$ and
            \begin{equation*}
                \overline{\mymathbb{f}}(\mathbb{Y})=\begin{bmatrix}
                    \overline{f}(y_0)\\
                    \vdots\\
                    \overline{f}(y_q)
                \end{bmatrix}.
            \end{equation*}
            Let $\overline{\nabla_Xf}(\mathbb{Y})=\overline{(A^\top)^\dagger}B\overline{\mymathbb{f}}(\mathbb{Y})$.  If $C\epsilon^*\kappa(A)<1$, then
            \begin{equation*}
                \left\|\nabla_X f(\mathbb{Y})-\overline{\nabla_X f}(\mathbb{Y})\right\| \le \sqrt{q+1}f_M\left\|A^\dagger\right\|\left\|B\right\|\epsilon^*\left(\frac{\sqrt{2}C\kappa(A)}{1-C\kappa(A)\epsilon^*} + \frac{C+1}{(1-C\epsilon^*)(1-C\kappa(A)\epsilon^*)}\right).
            \end{equation*}
        \end{theorem}
        \begin{proof}
            Since 
            \begin{equation*}
                \left\|\mymathbb{f}(\mathbb{Y})\right\| \le \sqrt{q+1}f_M,
            \end{equation*}
            Theorem \ref{thm:invne} gives
            \begin{align*}
                &~~\left\|(A^\top)^\dagger B\mymathbb{f}(\mathbb{Y})-\overline{(A^\top)^\dagger}B\mymathbb{f}(\mathbb{Y})\right\|\\
                &\le \sqrt{q+1}f_M\left\|A^\dagger\right\|\left\|B\right\|\epsilon^*\left(\frac{\sqrt{2}C\kappa(A)}{1-C\kappa(A)\epsilon^*}+\frac{C}{(1-C\epsilon^*)(1-C\kappa(A)\epsilon^*)}\right).
            \end{align*}
            Applying Theorem \ref{thm:fevalne} with $(A^\top)^\dagger=\overline{(A^\top)^\dagger}$ and Lemma \ref{lem:bdpinv}, we have
            \begin{equation*}
                \left\|\overline{(A^\top)^\dagger}B\mymathbb{f}(\mathbb{Y})-\overline{(A^\top)^\dagger}B\overline{\mymathbb{f}}(\mathbb{Y})\right\| \le \frac{\sqrt{q+1}f_M\|A^\dagger\|\|B\|\epsilon^*}{(1-C\epsilon^*)(1-C\kappa(A)\epsilon^*)}.
            \end{equation*}
            Using triangle inequality, we have
            \begin{align*}
                \left\|\nabla_X f(\mathbb{Y})-\overline{\nabla_X f}(\mathbb{Y})\right\| &= \left\|(A^\top)^\dagger B\mymathbb{f}(\mathbb{Y})-\overline{(A^\top)^\dagger}B\overline{\mymathbb{f}}(\mathbb{Y})\right\|\\
                &\le \left\|(A^\top)^\dagger B\mymathbb{f}(\mathbb{Y})-\overline{(A^\top)^\dagger}B\mymathbb{f}(\mathbb{Y})\right\| + \left\|\overline{(A^\top)^\dagger}B\mymathbb{f}(\mathbb{Y})-\overline{(A^\top)^\dagger}B\overline{\mymathbb{f}}(\mathbb{Y})\right\|,
            \end{align*}
            and the result follows.
    
        $\hfill\qed$    
        \end{proof}

        \begin{remark}
            As we mentioned in Subsection \ref{subsec:fevalfpe}, our analysis can be adapted to analyze the noise in objective functions instead of the floating point error in function evaluations.  If this is the case, then we note that the results in Theorem \ref{thm:invfevalne} may not hold directly.  This is because currently both Subsections \ref{subsec:invfpe} and \ref{subsec:fevalfpe} are using the machine precision $\epsilon^*$.  If the $\epsilon^*$ in Subsection \ref{subsec:fevalfpe} is replaced by the noise level, say $\epsilon_f$, then results similar to Theorem \ref{thm:invfevalne} can be established, which contain both $\epsilon^*$ and $\epsilon_f$.
        \end{remark}

\section{Examples}\label{sec:exps}
    As we mentioned at the beginning of Section \ref{sec:fpe}, all the GSG, GCSG, and GACSG satisfy the general form of simplex derivatives.  Therefore, the result of Theorem \ref{thm:invfevalne} can be directly applied to them to obtain error bounds of floating point errors in computation.  In this section, we provide details of how Theorem \ref{thm:invfevalne} can be applied to each of them.

    We note that Theorem \ref{thm:invfevalne} assumes that $A$ has full-column rank and $C\epsilon^*\kappa(A)<1$. The first assumption can be satisfied by picking sample points such that the columns of $A$ are linearly independent, which is reasonable since the columns of $A$ are normally the directions from $y_0$ to sample points or the linear combination of these directions.  The second assumption requires that the condition number of $A$ is not too large (approximately $10^{16}$), which is also reasonable since matrices with condition numbers around $10^{16}$ rarely occur in practice.
    
    \begin{example}\label{exp:GSGne}
        For the GSG with $\mathbb{Y}=\{y_0,y_0+d_1,\ldots,y_0+d_p\}$, we have $z=q=p$, $A=L=[d_1\cdots d_p]$, and $B=[-\mymathbb{1}~I_p]$.  Suppose that the assumptions of Theorem \ref{thm:invfevalne} hold.  
        
        Notice that 
        \begin{equation*}
            \left\|B\right\| = \sqrt{\lambda_{\max}\left(B^\top B\right)} = \sqrt{\lambda_{\max}\left(BB^\top\right)} = \sqrt{\lambda_{\max}\left(\mymathbb{1}\mymathbb{1}^\top + I_p\right)} = \sqrt{\lambda_{\max}\left(\mymathbb{1}\mymathbb{1}^\top\right) + 1} = \sqrt{p+1}.
        \end{equation*}
        Therefore, 
        \begin{align*}
            \left\|\nabla_S f(\mathbb{Y})-\overline{\nabla_S f}(\mathbb{Y})\right\| &\le \left(p+1\right)f_M\left\|L^\dagger\right\|\epsilon^*\left(\frac{\sqrt{2}C\kappa(L)}{1-C\kappa(L)\epsilon^*} + \frac{C+1}{(1-C\epsilon^*)(1-C\kappa(L)\epsilon^*)}\right)\\
            &= \left(p+1\right)f_M\left\|\widehat{L}^\dagger\right\|\epsilon^*\left(\frac{\sqrt{2}C\kappa(L)}{1-C\kappa(L)\epsilon^*} + \frac{C+1}{(1-C\epsilon^*)(1-C\kappa(L)\epsilon^*)}\right)\frac{1}{\Delta},
        \end{align*}
        where $\Delta=\overline{\mathrm{diam}}(\mathbb{Y})$ and $\widehat{L}=L\slash\Delta$.
    \end{example}
    
    \begin{example}\label{exp:GACSGne}
        For the GACSG with $\mathbb{Y}^+=\{y_0,y_0+d_1,\ldots,y_0+d_p\}$ and $\widetilde{\mathbb{Y}}=\{y_0-\widetilde{d}_1,\ldots,y_0-\widetilde{d}_p\}$, we have $\mathbb{Y}=\{y_0,y_0+d_1,\ldots,y_0+d_p,y_0-\widetilde{d}_1,\ldots,y_0-\widetilde{d}_p\}$, $z=p$, $q=2p$, $A=L^+D-\widetilde{L}=[k_1^2d_1+\widetilde{d}_1\cdots k_p^2d_p+\widetilde{d}_p]$, and \begin{equation*}
            B=\begin{bmatrix}
                1-k_1^2 & k_1^2 & \cdots & 0 & -1 & \cdots & 0\\
                \vdots & \vdots & \ddots & \vdots & \vdots & \ddots & \vdots\\
                1-k_p^2 & 0 & \cdots & k_p^2 & 0 & \cdots & -1
            \end{bmatrix}.
        \end{equation*}
        Suppose that the assumptions of Theorem \ref{thm:invfevalne} hold.  
        
        Let $\mymathbb{k}=[k_1^2\cdots k_p^2]^\top$.  Then $B=[\mymathbb{1}-\mymathbb{k}~D~-I_p]$.  Notice that
        \begin{align*}
            \left\|B\right\| &= \sqrt{\lambda_{\max}\left(BB^\top\right)}\\
            &= \sqrt{\lambda_{\max}\left(\left(\mymathbb{1}-\mymathbb{k}\right)\left(\mymathbb{1}-\mymathbb{k}\right)^\top + D^2 + I_p\right)}\\
            &\le \sqrt{\lambda_{\max}\left(\left(\mymathbb{1}-\mymathbb{k}\right)\left(\mymathbb{1}-\mymathbb{k}\right)^\top\right) + \max\left\{k_i^4+1\right\}}\\
            &= \sqrt{\sum\limits_{i=1}^p\left(1-k_i^2\right)^2 + \max\left\{k_i^4+1\right\}}
        \end{align*}
        Therefore,
        \begin{align*}
            &~~\left\|\nabla_{ACS} f(\mathbb{Y})-\overline{\nabla_{ACS} f}(\mathbb{Y})\right\|\\
            &\le \sqrt{2p+1}f_M\sqrt{\sum\limits_{i=1}^p\left(1-k_i^2\right)^2 + \max\left\{k_i^4+1\right\}}\left\|\left(L^+D-\widetilde{L}\right)^\dagger\right\|\epsilon^*\\
            &\hspace{2cm}\cdot\left(\frac{\sqrt{2}C\kappa(L^+D-\widetilde{L})}{1-C\kappa(L^+D-\widetilde{L})\epsilon^*} + \frac{C+1}{(1-C\epsilon^*)(1-C\kappa(L^+D-\widetilde{L})\epsilon^*)}\right)\\
            &= \sqrt{2p+1}f_M\sqrt{\sum\limits_{i=1}^p\left(1-k_i^2\right)^2 + \max\left\{k_i^4+1\right\}}\left\|\left(\widehat{L}^+D-\widehat{\widetilde{L}}\right)^\dagger\right\|\epsilon^*\\
            &\hspace{2cm}\cdot\left(\frac{\sqrt{2}C\kappa(L^+D-\widetilde{L})}{1-C\kappa(L^+D-\widetilde{L})\epsilon^*} + \frac{C+1}{(1-C\epsilon^*)(1-C\kappa(L^+D-\widetilde{L})\epsilon^*)}\right)\frac{1}{\Delta},
        \end{align*}
        where $\Delta=\overline{\mathrm{diam}}(\mathbb{Y}^+)$, $\widehat{L}^+=L^+\slash\Delta$ and $\widehat{\widetilde{L}}=\widetilde{L}\slash\Delta$.
    \end{example}

    \begin{example}\label{exp:GCSGne}
        For the GCSG with $\mathbb{Y}^+=\{y_0,y_0+d_1,\ldots,y_0+d_p\}$ and $\mathbb{Y}^-=\{y_0-d_1,\ldots,y_0-d_p\}$, we have $\mathbb{Y}=\{y_0,y_0+d_1,\ldots,y_0+d_p,y_0-d_1,\ldots,y_0-d_p\}$, $z=p$, $q=2p$, $A=2L^+=[2d_1\cdots 2d_p]$, and $B=[\mymathbb{0}~I_p~-I_p]$.  Suppose that the assumptions of Theorem \ref{thm:invfevalne} hold.  
        
        Taking all $k_i=1$, $D=I_p$, and $\widetilde{L}=-L^+$ in the Example \ref{exp:GACSGne}, we get
        \begin{align*}
            &~~\left\|\nabla_{CS} f(\mathbb{Y})-\overline{\nabla_{CS} f}(\mathbb{Y})\right\|\\
            &= \frac{\sqrt{4p+2}}{2}f_M\left\|\left(\widehat{L}^+\right)^\dagger\right\|\epsilon^*\left(\frac{\sqrt{2}C\kappa(L^+)}{1-C\kappa(L^+)\epsilon^*} + \frac{C+1}{(1-C\epsilon^*)(1-C\kappa(L^+)\epsilon^*)}\right)\frac{1}{\Delta},
        \end{align*}
        where $\Delta=\overline{\mathrm{diam}}(\mathbb{Y}^+)$ and $\widehat{L}^+=L^+\slash\Delta$.
    \end{example}

    We conclude with a final example, showing that if the structure of the simplex derivative is known in the first place, then a tighter error bound may exist.  This is because if matrix $B$ is known, then the structure of $B\mymathbb f(\mathbb{Y})$ can be exploited to cancel out some terms in the error analysis and improve the tightness of error bounds.
    \begin{example}[A tighter floating point error bound of the GSG]\label{exp:eb_GSGtighter}
        Consider the GSG with $\mathbb{Y}=\{y_0,y_0+d_1,\ldots,y_0+d_p\}$.  Suppose that the assumptions of Theorem \ref{thm:invfevalne} hold.  

        We denote $\epsilon_{y_0+d_i}$ by $\epsilon_i$, and $f(y_0+d_i)$ by $f_i$ for simplicity.  Notice that  
        \begin{align*}
            \left\|B\mymathbb{f}(\mathbb{Y})-B\overline{\mymathbb{f}}(\mathbb{Y})\right\| &= \sqrt{\sum\limits_{i=1}^p\left|\overline{f}(y_0+d_i)-\overline{f}(y_0)-\left(f(y_0+d_i)-f(y_0)\right)\right|^2}\\
            &= \sqrt{\sum\limits_{i=1}^p\left|\epsilon_if_i-\epsilon_0f_0\right|^2}\\
            &\le \sqrt{\sum\limits_{i=1}^p\left(\epsilon_i^2f_i^2+\epsilon_0^2f_0^2+2\left|\epsilon_if_i\epsilon_0f_0\right|\right)}\\
            &\le 2\sqrt{p}f_M\epsilon^*
        \end{align*}
        and
        \begin{equation*}
            \left\|B\mymathbb{f}(\mathbb{Y})\right\| = \sqrt{\sum\limits_{i=1}^p\left|f(y_0+d_i)-f(y_0)\right|^2} \le \sqrt{\sum\limits_{i=1}^p4f_M^2} = 2\sqrt{p}f_M.
        \end{equation*}
        Using triangle inequality and Lemmas \ref{lem:invne} and \ref{lem:bdpinv}, we have
        \begin{align*}
            &~~\left\|\nabla_S f(\mathbb{Y})-\overline{\nabla_S f}(\mathbb{Y})\right\|\\
            &\le \left\|(L^\top)^\dagger B\mymathbb{f}(\mathbb{Y})-\overline{(L^\top)^\dagger}B\mymathbb{f}(\mathbb{Y})\right\| + \left\|\overline{(L^\top)^\dagger}B\mymathbb{f}(\mathbb{Y})-\overline{(L^\top)^\dagger}B\overline{\mymathbb{f}}(\mathbb{Y})\right\|\\
            &\le \left\|(L^\top)^\dagger - \overline{(L^\top)^\dagger}\right\|\left\|B\mymathbb{f}(\mathbb{Y})\right\| + \left\|\overline{(L^\top)^\dagger}\right\|\left\|B\mymathbb{f}(\mathbb{Y})-B\overline{\mymathbb{f}}(\mathbb{Y})\right\|\\
            &\le 2\sqrt{p}f_M\left\|L^\dagger\right\|\epsilon^*\left(\frac{\sqrt{2}C\kappa(L)}{1-C\kappa(L)\epsilon^*}+\frac{C}{(1-C\epsilon^*)(1-C\kappa(L)\epsilon^*)}\right) + \frac{2\sqrt{p}f_M\left\|L^\dagger\right\|\epsilon^*}{(1-C\epsilon^*)(1-C\kappa(L)\epsilon^*)}\\
            &= 2\sqrt{p}f_M\left\|L^\dagger\right\|\epsilon^*\left(\frac{\sqrt{2}C\kappa(L)}{1-C\kappa(L)\epsilon^*} + \frac{C+1}{\left(1-C\epsilon^*\right)\left(1-C\kappa(L)\epsilon^*\right)}\right)\\
            &= 2\sqrt{p}f_M\left\|\widehat{L}^\dagger\right\|\epsilon^*\left(\frac{\sqrt{2}C\kappa(L)}{1-C\kappa(L)\epsilon^*} + \frac{C+1}{\left(1-C\epsilon^*\right)\left(1-C\kappa(L)\epsilon^*\right)}\right)\frac{1}{\Delta}.
        \end{align*}
    \end{example}
    \begin{remark}
        To compare the error bound of Example \ref{exp:eb_GSGtighter} to that of Example \ref{exp:GSGne}, we denote
        \begin{equation*}
            K = f_M\left\|\widehat{L}^\dagger\right\|\epsilon^*\left(\frac{\sqrt{2}C\kappa(L)}{1-C\kappa(L)\epsilon^*} + \frac{C+1}{\left(1-C\epsilon^*\right)\left(1-C\kappa(L)\epsilon^*\right)}\right)\frac{1}{\Delta}.
        \end{equation*}
        Then, Example \ref{exp:eb_GSGtighter} gives
        \begin{align*}
            \left\|\nabla_S f(\mathbb{Y})-\overline{\nabla_S f}(\mathbb{Y})\right\| \le 2\sqrt{p}K
        \end{align*}
        and Example \ref{exp:GSGne} gives
        \begin{equation*}
            \left\|\nabla_S f(\mathbb{Y})-\overline{\nabla_S f}(\mathbb{Y})\right\| \le \left(p+1\right)K.
        \end{equation*}
        We see that the error bound given by Example \ref{exp:eb_GSGtighter} is strictly tighter for all $p\ge 2$.  This implies that tighter error bounds may exist if we specifically analyze one type of simplex gradient. 
        % For the GACSG, the error bound above
        % \begin{align*}
        %     & \sqrt{2}f_M\sqrt{\sum\limits_{i=1}^pk_i^4+\sum\limits_{i=1}^p\left|k_i^4-1\right|+p}\left\|\left(L^+D-\widetilde{L}\right)^\dagger\right\|\epsilon^*\\
        %     &~~~~\cdot\left(\frac{\sqrt{2}C\kappa_2}{1-C\kappa_2\epsilon^*} + \frac{C+1}{(1-C\epsilon^*)(1-C\kappa_2\epsilon^*)} - \frac{C\kappa_2\epsilon^*}{1-C\kappa_2\epsilon^*} - \frac{C\epsilon^*}{(1-C\epsilon^*)(1-C\kappa_2\epsilon^*)}\right),
        % \end{align*}
        % and from the general form
        % \begin{equation*}
        %     \sqrt{2p+1}f_M\sqrt{\sum\limits_{i=1}^p\left(1-k_i^2\right)^2 + \max\left\{k_i^4+1\right\}}\left\|\left(L^+D-\widetilde{L}\right)^\dagger\right\|\epsilon^*\left(\frac{\sqrt{2}C\kappa_2}{1-C\kappa_2\epsilon^*} + \frac{C+1}{\left(1-C\epsilon^*\right)\left(1-C\kappa_2\epsilon^*\right)}\right).
        % \end{equation*}
    \end{remark}

\section{General error bound with floating point errors and minimal choice of $\Delta$}\label{sec:deltamin}
    The approximation error bounds provided in Section \ref{sec:def&appeb} suggest that the accuracy of approximation continuously improves as the sample set approximate diameter $\Delta$ goes to zero.  However, the floating point error bounds in Section \ref{sec:fpe} show that the floating point error increases as $\Delta$ goes to zero.  This trade-off implies that, as $\Delta$ goes to zero, there should be a turning point, i.e., a minimal choice of $\Delta$, after which floating point errors will take over and cause a decrease in approximation accuracy.

    In this section, we first combine the results of Sections \ref{sec:def&appeb} and \ref{sec:fpe} and provide a general error bound between the true gradient and the computed simplex derivatives.  Then, we examine the minimal choice of $\Delta$ for simplex derivatives.  As examples, we specialize our results to the GSG, GCSG, and GACSG to show how our results can be applied to them.

    Notice that the approximation errors in Section \ref{sec:def&appeb} only provide accuracy guarantees on the column spaces. Therefore, we can only obtain an error estimate between the true gradient and the computed simplex derivatives on these column spaces, which is shown in the following lemma.
    \begin{lemma}\label{lem:eb_general}
        Let $\overline{\nabla_Xf}(\mathbb{Y})=\overline{(A^\top)^\dagger}B\overline{\mymathbb{f}}(\mathbb{Y})$.  Suppose that
        \begin{equation*}
            \left\|\mathrm{proj}_{\mathrm{col}(A)}(\nabla f(y_0))-\mathrm{proj}_{\mathrm{col}(A)}(\nabla_X f(\mathbb{Y}))\right\| \le \kappa^{\mathrm{ae}}(y_0)\Delta^{N^{\mathrm{ae}}}
        \end{equation*}
        and
        \begin{equation*}
            \left\|\nabla_X f(\mathbb{Y})-\overline{\nabla_X f}(\mathbb{Y})\right\| \le \kappa^{\mathrm{fpe}}(y_0)\Delta^{-N^{\mathrm{fpe}}},
        \end{equation*}
        where $\Delta=\overline{\mathrm{diam}}(\mathbb{Y})$, and $\kappa^{\mathrm{ae}}(y_0)>0,\kappa^{\mathrm{fpe}}(y_0)>0,N^{\mathrm{ae}}\ge 1,N^{\mathrm{fpe}}> 0$ are independent of $\Delta$.  Then we have
        \begin{equation*}
            \left\|\mathrm{proj}_{\mathrm{col}(A)}(\nabla f(y_0))-\mathrm{proj}_{\mathrm{col}(A)}(\overline{\nabla_X f}(\mathbb{Y}))\right\| \le \kappa^{\mathrm{ae}}(y_0)\Delta^{N^{\mathrm{ae}}} + \kappa^{\mathrm{fpe}}(y_0)\Delta^{-N^\mathrm{fpe}}.
        \end{equation*}
    \end{lemma}
    \begin{proof}
        Notice that
        \begin{align*}
            \left\|\mathrm{proj}_{\mathrm{col}(A)}(\nabla_X f(\mathbb{Y}))-\mathrm{proj}_{\mathrm{col}(A)}(\overline{\nabla_X f}(\mathbb{Y}))\right\| &= \left\|\mathrm{proj}_{\mathrm{col}(A)}(\nabla_X f(\mathbb{Y})-\overline{\nabla_X f}(\mathbb{Y}))\right\|\\
            &\le \left\|\nabla_X f(\mathbb{Y})-\overline{\nabla_X f}(\mathbb{Y})\right\|\\
            &\le \kappa^{\mathrm{fpe}}(y_0)\Delta^{-N^{\mathrm{fpe}}}.
        \end{align*}
        The result follows from the triangle inequality.

    $\hfill\qed$   
    \end{proof}
    
    In particular, combining the results in Sections \ref{sec:def&appeb} to \ref{sec:fpe}, we have the GSG and GCSG satisfy the assumptions of Lemma \ref{lem:eb_general}.  Based on the general error bound between the true gradient and the computed simplex derivatives, we provide the minimal choice of $\Delta$ as follows.  
    \begin{theorem}[Minimal choice of $\Delta$]\label{thm:deltamin}
        Let $\overline{\nabla_Xf}(\mathbb{Y})=\overline{(A^\top)^\dagger}B\overline{\mymathbb{f}}(\mathbb{Y})$.  Suppose that 
        \begin{equation*}
            \left\|\mathrm{proj}_{\mathrm{col}(A)}(\nabla f(y_0))-\mathrm{proj}_{\mathrm{col}(A)}(\overline{\nabla_X f}(\mathbb{Y}))\right\| \le \kappa^{\mathrm{ae}}(y_0)\Delta^{N^{\mathrm{ae}}} + \kappa^{\mathrm{fpe}}(y_0)\Delta^{-N^\mathrm{fpe}},
        \end{equation*}
        where $\Delta=\overline{\mathrm{diam}}(\mathbb{Y})$, and $\kappa^{\mathrm{ae}}(y_0)>0,\kappa^{\mathrm{fpe}}(y_0)>0,N^{\mathrm{ae}}\ge 1,N^{\mathrm{fpe}}> 0$ are independent of $\Delta$.  Then the minimal choice of $\Delta=\overline{\mathrm{diam}}(\mathbb{Y})$ is
        \begin{equation*}
            \Delta_X^{\mathrm{min}} = \left(\frac{N^{\mathrm{fpe}}\kappa^{\mathrm{fpe}}(y_0)}{N^{\mathrm{ae}}\kappa^{\mathrm{ae}}(y_0)}\right)^{1\slash(N^\mathrm{ae}+N^\mathrm{fpe})}.
        \end{equation*}
    \end{theorem}
    \begin{proof}
        We denote 
        \begin{equation*}
            F(\Delta)=\kappa^{\mathrm{ae}}(y_0)\Delta^{N^{\mathrm{ae}}} + \kappa^{\mathrm{fpe}}(y_0)\Delta^{-N^\mathrm{fpe}}~~\text{and}~~\mathcal{D}=\{\Delta:\Delta>0\}.
        \end{equation*}
        Then 
        \begin{equation*}
            \Delta_X^{\mathrm{min}} \in \argmin_{\mathcal{D}}F(\Delta).
        \end{equation*}
        Since $\kappa^{\mathrm{ae}}(y_0)>0,\kappa^{\mathrm{fpe}}(y_0)>0,N^{\mathrm{ae}}\ge 1,N^{\mathrm{fpe}}> 0$, we have $\kappa^{\mathrm{ae}}(y_0)\Delta^{N^{\mathrm{ae}}}$ is convex on $\mathcal{D}$ and $\kappa^{\mathrm{fpe}}(y_0)\Delta^{-N^\mathrm{fpe}}$ is strictly convex on $\mathcal{D}$.  Therefore, $\Delta_X^{\mathrm{min}}$ is the unique minimizer of $F(\Delta)$ on $\mathcal{D}$ and satisfies
        \begin{equation*}
            \frac{\mathrm{d}F}{\mathrm{d}\Delta}(\Delta_X^{\mathrm{min}}) = N^{\mathrm{ae}}\kappa^{\mathrm{ae}}(y_0)(\Delta_X^{\mathrm{min}})^{N^{\mathrm{ae}}-1} - N^\mathrm{fpe}\kappa^{\mathrm{fpe}}(y_0)(\Delta_X^{\mathrm{min}})^{-N^\mathrm{fpe}-1} = 0.
        \end{equation*}
        Solving the equation above for $\Delta_X^{\mathrm{min}}$, we obtain
        \begin{equation*}
            \Delta_X^{\mathrm{min}} = \left(\frac{N^{\mathrm{fpe}}\kappa^{\mathrm{fpe}}(y_0)}{N^{\mathrm{ae}}\kappa^{\mathrm{ae}}(y_0)}\right)^{1\slash(N^\mathrm{ae}+N^\mathrm{fpe})}.
        \end{equation*}

    $\hfill\qed$   
    \end{proof}

    Now we present examples applying the results of Lemma \ref{lem:eb_general} and Theorem \ref{thm:deltamin} to the GSG and GCSG to obtain the error estimate and the minimal choice of $\Delta$.
    \begin{example}\label{exp:GSGdelta}
        Consider the GSG with $\mathbb{Y}=\{y_0,y_0+d_1,\ldots,y_0+d_p\}$.  Suppose that the assumptions of Theorems \ref{thm:GSGeb} and \ref{thm:invfevalne} hold.  Then
        \begin{align*}
            &~~\left\|\mathrm{proj}_{\mathrm{col}(L)}(\nabla f(y_0))-\mathrm{proj}_{\mathrm{col}(L)}(\nabla_{S}f(\mathbb{Y}))\right\|\\
            &\le \frac{\nu\sqrt{p}}{2}\left\|\widehat{L}^\dagger\right\|\Delta + \left(p+1\right)f_M\left\|\widehat{L}^\dagger\right\|\epsilon^*\left(\frac{\sqrt{2}C\kappa(L)}{1-C\kappa(L)\epsilon^*} + \frac{C+1}{(1-C\epsilon^*)(1-C\kappa(L)\epsilon^*)}\right)\frac{1}{\Delta},
        \end{align*}
        where $\Delta=\overline{\mathrm{diam}}(\mathbb{Y})$ and $\widehat{L}=L\slash\Delta$.  Hence, 
        \begin{align*}
            \kappa^{\mathrm{ae}}(y_0) &= \frac{\nu\sqrt{p}}{2}\left\|\widehat{L}^\dagger\right\|,\\
            \kappa^{\mathrm{fpe}}(y_0) &= \left(p+1\right)f_M\left\|\widehat{L}^\dagger\right\|\epsilon^*\left(\frac{\sqrt{2}C\kappa(L)}{1-C\kappa(L)\epsilon^*} + \frac{C+1}{(1-C\epsilon^*)(1-C\kappa(L)\epsilon^*)}\right),\\
            N^{\mathrm{ae}} &= 1,\\
            N^{\mathrm{fpe}} &= 1,
        \end{align*}
        which gives
        \begin{equation*}
            \Delta_S^{\mathrm{min}} = \left(\frac{2\left(p+1\right)f_M\epsilon^*}{\nu\sqrt{p}}\left(\frac{\sqrt{2}C\kappa(L)}{1-C\kappa(L)\epsilon^*} + \frac{C+1}{(1-C\epsilon^*)(1-C\kappa(L)\epsilon^*)}\right)\right)^{1\slash 2}.
        \end{equation*}

        If we use the tighter floating error bound in Example \ref{exp:eb_GSGtighter}, then we have
        \begin{align*}
            &~~\left\|\mathrm{proj}_{\mathrm{col}(L)}(\nabla f(y_0))-\mathrm{proj}_{\mathrm{col}(L)}(\nabla_{S}f(\mathbb{Y}))\right\|\\
            &\le \frac{\nu\sqrt{p}}{2}\left\|\widehat{L}^\dagger\right\|\Delta + 2\sqrt{p}f_M\left\|\widehat{L}^\dagger\right\|\epsilon^*\left(\frac{\sqrt{2}C\kappa(L)}{1-C\kappa(L)\epsilon^*} + \frac{C+1}{\left(1-C\epsilon^*\right)\left(1-C\kappa(L)\epsilon^*\right)}\right)\frac{1}{\Delta},
        \end{align*}
        where $\Delta=\overline{\mathrm{diam}}(\mathbb{Y})$ and $\widehat{L}=L\slash\Delta$.  Hence, $\kappa^{\mathrm{ae}}(y_0)$, $N^{\mathrm{ae}}$, and $N^{\mathrm{fpe}}$ are the same as above, and
        \begin{equation*}
            \kappa^{\mathrm{fpe}}(y_0) = 2\sqrt{p}f_M\left\|\widehat{L}^\dagger\right\|\epsilon^*\left(\frac{\sqrt{2}C\kappa(L)}{1-C\kappa(L)\epsilon^*} + \frac{C+1}{\left(1-C\epsilon^*\right)\left(1-C\kappa(L)\epsilon^*\right)}\right),
        \end{equation*}
        which gives
        \begin{equation*}
            \Delta_S^{\mathrm{min}} = \left(\frac{4f_M\epsilon^*}{\nu}\left(\frac{\sqrt{2}C\kappa(L)}{1-C\kappa(L)\epsilon^*} + \frac{C+1}{(1-C\epsilon^*)(1-C\kappa(L)\epsilon^*)}\right)\right)^{1\slash 2}.
        \end{equation*}
    
        For a rough estimate of each of the above $\Delta_S^{\mathrm{min}}$, we suppose that $\nu=f_M=C=1$.  Then, from Example \ref{exp:GSGne} we have
        \begin{equation*}
            \Delta_S^{\mathrm{min}} = \left(\frac{2\left(p+1\right)\epsilon^*}{\sqrt{p}}\left(\frac{\sqrt{2}\kappa(L)}{1-\kappa(L)\epsilon^*} + \frac{2}{(1-\epsilon^*)(1-\kappa(L)\epsilon^*)}\right)\right)^{1\slash 2}.
        \end{equation*}
        If we use the tighter floating error bound in Example \ref{exp:eb_GSGtighter}, then we have
        \begin{equation*}
            \Delta_S^{\mathrm{min}} = \left(4\epsilon^*\left(\frac{\sqrt{2}\kappa(L)}{1-\kappa(L)\epsilon^*} + \frac{2}{(1-\epsilon^*)(1-\kappa(L)\epsilon^*)}\right)\right)^{1\slash 2}.
        \end{equation*}
    \end{example}

    \begin{example}\label{exp:GCSGdelta}
        Consider the GCSG with $\mathbb{Y}^+=\{y_0,y_0+d_1,\ldots,y_0+d_p\},\mathbb{Y}^-=\{y_0-d_1,\ldots,y_0-d_p\}$, and $\mathbb{Y}=\{y_0,y_0+d_1,\ldots,y_0+d_p,y_0-d_1,\ldots,y_0-d_p\}$.  Suppose that the assumptions of Theorems~\ref{thm:GCSGeb} and \ref{thm:invfevalne} hold.  Then
        \begin{align*}
            &~~\left\|\mathrm{proj}_{\mathrm{col}(L^+)}(\nabla f(y_0))-\mathrm{proj}_{\mathrm{col}(L^+)}(\nabla_{CS}f(\mathbb{Y}))\right\|\\
            &\le \frac{\nu\sqrt{p}}{6}\left\|\left(\widehat{L}^+\right)^\dagger\right\|\Delta^2 + \frac{\sqrt{4p+2}}{2}f_M\left\|\left(\widehat{L}^+\right)^\dagger\right\|\epsilon^*\left(\frac{\sqrt{2}C\kappa(L^+)}{1-C\kappa(L^+)\epsilon^*} + \frac{C+1}{(1-C\epsilon^*)(1-C\kappa(L^+)\epsilon^*)}\right)\frac{1}{\Delta},
        \end{align*}
        where $\Delta=\overline{\mathrm{diam}}(\mathbb{Y}^+)$ and $\widehat{L}^+=L^+\slash\Delta$.  Hence, 
        \begin{align*}
            \kappa^{\mathrm{ae}}(y_0) &= \frac{\nu\sqrt{p}}{6}\left\|\left(\widehat{L}^+\right)^\dagger\right\|,\\
            \kappa^{\mathrm{fpe}}(y_0) &= \frac{\sqrt{4p+2}}{2}f_M\left\|\left(\widehat{L}^+\right)^\dagger\right\|\epsilon^*\left(\frac{\sqrt{2}C\kappa(L^+)}{1-C\kappa(L^+)\epsilon^*} + \frac{C+1}{(1-C\epsilon^*)(1-C\kappa(L^+)\epsilon^*)}\right),\\
            N^{\mathrm{ae}} &= 2,\\
            N^{\mathrm{fpe}} &= 1,
        \end{align*}
        which gives
        \begin{equation*}
            \Delta_{CS}^{\mathrm{min}} = \left(\frac{3\sqrt{4p+2}f_M\epsilon^*}{\nu\sqrt{p}}\left(\frac{\sqrt{2}C\kappa(L^+)}{1-C\kappa(L^+)\epsilon^*} + \frac{C+1}{(1-C\epsilon^*)(1-C\kappa(L^+)\epsilon^*)}\right)\right)^{1\slash 3}.
        \end{equation*}

        For a rough estimate of $\Delta_{CS}^{\mathrm{min}}$, we suppose that $\nu=f_M=C=1$.  Then
        \begin{equation*}
            \Delta_{CS}^{\mathrm{min}} = \left(\frac{3\sqrt{4p+2}\epsilon^*}{\sqrt{p}}\left(\frac{\sqrt{2}\kappa(L^+)}{1-\kappa(L^+)\epsilon^*} + \frac{2}{(1-\epsilon^*)(1-\kappa(L^+)\epsilon^*)}\right)\right)^{1\slash 3}.
        \end{equation*}
    \end{example}

    For the GACSG, the error bound is not in the general form required by Theorem \ref{thm:deltamin}, as its error bound in Section \ref{sec:def&appeb} has two terms with different order of $\Delta$.  However, as shown in the following example, the error bound of the GACSG can be shaped into the general form under reasonable assumptions.
    \begin{example}\label{exp:GACSGdelta}
        Consider the GACSG with $\mathbb{Y}^+=\{y_0,y_0+d_1,\ldots,y_0+d_p\},\widetilde{\mathbb{Y}}=\{y_0-\widetilde{d}_1,\ldots,y_0-\widetilde{d}_p\}$, and $\mathbb{Y}=\{y_0,y_0+d_1,\ldots,y_0+d_p,y_0-\widetilde{d}_1,\ldots,y_0-\widetilde{d}_p\}$.  Suppose that the assumptions of Theorems~\ref{thm:GACSGeb} and \ref{thm:invfevalne} hold.  Then
        \begin{equation*}
            \left\|\mathrm{proj}_{\mathrm{col}(L^+D-\widetilde{L})}(\nabla f(y_0))-\mathrm{proj}_{\mathrm{col}(L^+D-\widetilde{L})}(\nabla_{ACS}f(\mathbb{Y}))\right\| \le \kappa_1^{\mathrm{ae}}(y_0)\Delta + \kappa_2^{\mathrm{ae}}(y_0)\Delta^{2} + \kappa^{\mathrm{fpe}}(y_0)\Delta^{-1},
        \end{equation*}
        where 
        \begin{align*}
            \kappa_1^{\mathrm{ae}}(y_0) &= \frac{1}{2}\max\left\{k_i^2\right\}\sqrt{p}\left\|\left(\widehat{L}^+D-\widehat{\widetilde{L}}\right)^\dagger\right\|\left(K_1\max\left\{\sin\theta_i\right\} + K_2\max\left\{\sqrt{1-\cos\theta_i}\right\}\right),\\
            \kappa_2^{\mathrm{ae}}(y_0) &= \frac{\nu}{6}\max\left\{k_i^2\left(1+k_i\right)\right\}\sqrt{p}\left\|\left(\widehat{L}^+D-\widehat{\widetilde{L}}\right)^\dagger\right\|\\
            \kappa^{\mathrm{fpe}}(y_0) &=\sqrt{2p+1}f_M\sqrt{\sum\limits_{i=1}^p\left(1-k_i^2\right)^2 + \max\left\{k_i^4+1\right\}}\left\|\left(\widehat{L}^+D-\widehat{\widetilde{L}}\right)^\dagger\right\|\epsilon^*\\
            &\hspace{2cm}\cdot\left(\frac{\sqrt{2}C\kappa(L^+D-\widetilde{L})}{1-C\kappa(L^+D-\widetilde{L})\epsilon^*} + \frac{C+1}{(1-C\epsilon^*)(1-C\kappa(L^+D-\widetilde{L})\epsilon^*)}\right),
        \end{align*}
        where $\Delta=\overline{\mathrm{diam}}(\mathbb{Y}^+)$, $\widehat{L}^+=L^+\slash\Delta$ and $\widehat{\widetilde{L}}=\widetilde{L}\slash\Delta$.

        If all $\theta_i$ are small enough such that $\kappa_1^{\mathrm{ae}}(y_0) << \kappa_2^{\mathrm{ae}}(y_0)\Delta$, then we can apply Theorem \ref{thm:deltamin} with $\kappa^{\mathrm{ae}}(y_0)=2\kappa_2^{\mathrm{ae}}(y_0), N^{\mathrm{ae}}=2, N^{\mathrm{fpe}}=1$ to get
        \begin{align*}
            \Delta_{ACS}^{\mathrm{min}} &= \left(\frac{3\sqrt{2p+1}f_M\sqrt{\sum\limits_{i=1}^p\left(1-k_i^2\right)^2 + \max\left\{k_i^4+1\right\}}\epsilon^*}{\nu\max\left\{k_i^2\left(1+k_i\right)\right\}\sqrt{p}}\right.\\
            &\hspace{2cm}\left.\cdot\left(\frac{\sqrt{2}C\kappa(L^+D-\widetilde{L})}{1-C\kappa(L^+D-\widetilde{L})\epsilon^*} + \frac{C+1}{(1-C\epsilon^*)(1-C\kappa(L^+D-\widetilde{L})\epsilon^*)}\right)\right)^{1\slash 3}.
        \end{align*}

        Conversely, if $\Delta$ is small enough such that $\kappa_2^{\mathrm{ae}}(y_0)\Delta << \kappa_1^{\mathrm{ae}}(y_0)$, then
        we can apply Theorem~\ref{thm:deltamin} with $\kappa^{\mathrm{ae}}(y_0)=2\kappa_1^{\mathrm{ae}}(y_0), N^{\mathrm{ae}}=1, N^{\mathrm{fpe}}=1$ to get
        \begin{align*}
            \Delta_{ACS}^{\mathrm{min}} &= \left(\frac{\sqrt{2p+1}f_M\sqrt{\sum\limits_{i=1}^p\left(1-k_i^2\right)^2 + \max\left\{k_i^4+1\right\}}\epsilon^*}{\max\left\{k_i^2\right\}\sqrt{p}}\right.\\
            &\hspace{2cm}\left.\cdot\left(\frac{\sqrt{2}C\kappa(L^+D-\widetilde{L})}{1-C\kappa(L^+D-\widetilde{L})\epsilon^*} + \frac{C+1}{(1-C\epsilon^*)(1-C\kappa(L^+D-\widetilde{L})\epsilon^*)}\right)\right)^{1\slash 2}.
        \end{align*}
        
        For a rough estimate of $\Delta_{ACS}^{\mathrm{min}}$, we suppose that all $k_i=1$ and $\nu=f_M=C=1$.  If $\kappa_1^{\mathrm{ae}}(y_0) << \kappa_2^{\mathrm{ae}}(y_0)\Delta$, then
        \begin{equation*}
            \Delta_{ACS}^{\mathrm{min}} = \left(\frac{3\sqrt{4p+2}\epsilon^*}{2\sqrt{p}}\left(\frac{\sqrt{2}\kappa(L^+D-\widetilde{L})}{1-\kappa(L^+D-\widetilde{L})\epsilon^*} + \frac{2}{(1-\epsilon^*)(1-\kappa(L^+D-\widetilde{L})\epsilon^*)}\right)\right)^{1\slash 3}.
        \end{equation*}
        If $\kappa_2^{\mathrm{ae}}(y_0)\Delta << \kappa_1^{\mathrm{ae}}(y_0)$, then
        \begin{equation*}
            \Delta_{ACS}^{\mathrm{min}} = \left(\frac{\sqrt{4p+2}\epsilon^*}{\sqrt{p}}\left(\frac{\sqrt{2}\kappa(L^+D-\widetilde{L})}{1-\kappa(L^+D-\widetilde{L})\epsilon^*} + \frac{2}{(1-\epsilon^*)(1-\kappa(L^+D-\widetilde{L})\epsilon^*)}\right)\right)^{1\slash 2}.
        \end{equation*}
    \end{example}

\section{Conclusion}\label{sec:concl}
    In this paper, we reviewed some popular gradient approximation techniques used in derivative-free optimization including the GSG and GCSG.  The main focus of this paper is to present a general framework for floating point error analysis of simplex derivatives, which does not depend on the choice of the simplex derivative as long as they are in a general form.  Our framework explains how the two major sources of floating point errors, i.e., pseudo-inverse computation and objective function evaluation, influence the approximation accuracy of simplex derivatives. Notice that our analysis in Section \ref{sec:fpe} can be easily adapted to analyze the influence of noise in the objective function.  In addition, we also advanced and analyzed GACSG, which generalizes the adapted centred simplex gradient proposed in \cite{chen2023adapting} to allow the sample set to contain any number of points.  We applied our framework to the GSG, GCSG, and GACSG and gave suggestions on the minimal choice of $\Delta$, the approximate diameter of the sample set.

    We note that error bounds tighter than expected from our framework may exist.  As shown in Example \ref{exp:eb_GSGtighter}, if we analyze the floating point error of a specific type of simplex derivative, then the structure of this simplex derivative can be exploited to get a tighter upper bound.  Future research may consider the general structure of a class of simplex derivatives and improve our framework for that class.
    \bigskip

\noindent{\bf Disclosure Statement:} The authors have no competing interests. 

\normalsize
\bibliographystyle{siam}
\bibliography{references}
\end{document}